\documentclass[letterpaper, 10 pt, conference]{ieeeconf}  

\IEEEoverridecommandlockouts                              
\newcommand{\bC}{\mathbb{C}}

\usepackage{graphics} 
\usepackage{amsmath} 
\usepackage{amssymb}  

\usepackage{mathtools}
\usepackage{cite,color}
\usepackage[american,fulldiode]{circuitikz} 

\usepackage{etoolbox}

\AtEndEnvironment{thebibliography}{

\bibitem{SottileBook}
F. Sottile, \emph{Real Solutions to Equations from Geometry},
University Lecture Series, ~57, American Mathematical Society, Providence, RI, 2011.

\bibitem{Harris1979}
J. Harris, ``{Galois groups of enumerative problems},''
\emph{Duke Math. J.}, vol.~46, no.~4, pp.
  685--724, 1979.

\bibitem{BertiniBook}
D.J. Bates, J.D. Hauenstein, A.J. Sommese, C.W. Wampler,
\emph{Numerically Solving Polynomial Systems with Bertini},
Software, Environments, and Tools, 25,
Society for Industrial and Applied Mathematics (SIAM), Philadelphia, PA, 2013.

\bibitem{LeykinSottile}
A. Leykin and F. Sottile, ``{Galois groups of {S}chubert problems via homotopy computation},''
\emph{Math. Comp.}, vol.~78, no.~267, 1749--1765, 2009.

\bibitem{GaloisHRS}
J.D. Hauenstein, J.I. Rodriguez, and F. Sottile,
``{Numerical computation of Galois groups},''
In prepartion, 2016.
}

\title{\LARGE \bf
Investigating the Maximum Number of Real Solutions to the \\ Power Flow Equations: 
Analysis of Lossless Four-Bus Systems
}

\author{Daniel K. Molzahn$^{1}$, Matthew Niemerg$^{2}$, Dhagash Mehta$^{3}$, and Jonathan D. Hauenstein$^{4}$
\thanks{$^{1}$Argonne National Laboratory, Energy Systems Division.
        {\tt\small dmolzahn@anl.gov}}%
\thanks{$^{2}$IBM T.J. Watson Research Center.  Supported
in part by the Fields Institute.  {\tt\small mniemer@us.ibm.com}}%
\thanks{$^{3}$University of Notre Dame, Dept. of Applied and Computational Mathematics and Statistics. 
Supported in part by NSF ECCS 1509036.
        {\tt\small dmehta@nd.edu}}%
\thanks{$^{4}$University of Notre Dame, Dept. of Applied and Computational Mathematics and Statistics. 
Supported in part by NSF ACI 1460032, Sloan Research Fellowship, and Army Young Investigator Program (YIP).
        {\tt\small hauenstein@nd.edu}}%
}

\begin{document}

\maketitle
\thispagestyle{empty}
\pagestyle{empty}

\begin{abstract}
The power flow equations model the steady-state relationship between the power injections and voltage phasors in an electric power system. By separating the real and imaginary components of the voltage phasors, the power flow equations can be formulated as a system of quadratic polynomials. Only the real solutions to these polynomial equations are physically meaningful. This paper focuses on the maximum number of real solutions to the power flow equations. An upper bound on the number of real power flow solutions commonly used in the literature is the maximum number of complex solutions. There exist two- and three-bus systems for which all complex solutions are real. It is an open question whether this is also the case for larger systems. This paper investigates four-bus systems using techniques from numerical algebraic geometry and conjectures a negative answer to this question. In particular, this paper studies lossless, four-bus systems composed of PV buses connected by lines with arbitrary susceptances. Computing the Galois group, which is degenerate, enables conversion of the problem of counting the number of real solutions to the power flow equations into counting the number of positive roots of a univariate sextic polynomial. From this analysis, it is conjectured that the system has at most 16 real solutions, which is strictly less than the maximum number of complex solutions, namely 20.  We also provide explicit parameter values where this system has 16 real solutions so that the conjectured upper bound is achievable.
\end{abstract}

\section{Introduction}

The power flow equations are at the heart of many electric power system computations. These equations model the steady-state relationship between the power injections and voltage phasors in a power system. Power flow solutions correspond to the equilibria of the differential-algebraic equations describing the power system dynamics. The power flow equations also form the key constraints in many power system optimization problems, such as optimal power flow~\cite{florin2016}, unit commitment~\cite{tahanan2015}, and state estimation~\cite{huang2012}.

Despite their importance to power system analyses, there are open theoretical questions regarding power flow solution characteristics. This paper focuses on an open question regarding the maximum number of power flow solutions.

For typical operating conditions, there is a single ``high-voltage'' power flow solution which corresponds to a stable equilibrium of reasonable dynamical models. There exists a voluminous literature of mature algorithms for calculating this solution for large-scale systems (e.g., Newton-Raphson~\cite{tinney1967} and Gauss-Seidel~\cite{glimn1957}). The performance of Newton-based algorithms is sensitive to the initialization. While reasonable initializations (e.g., the solution to a related problem or a ``flat start'' of $1\angle 0^\circ$ per unit voltages) typically converge to the high-voltage solution, Newton-based methods have fractal 
basins of attraction~\cite{thorp1989}. Initializing Newton-based methods is thus challenging when parameters move outside typical operating ranges. Accordingly, 
significant research has been focused on modifications and alternatives to Newton-based methods. See~\cite{acc2016_tutorial} for a discussion of some recent developments.

It is well known that the power flow equations may have multiple solutions~\cite{klos1975}. While the high-voltage solution is often of primary interest, other solutions are relevant for certain stability assessments~\cite{venikov1975,tamura1983,ribbens-pavella1985,chiang1987,overbye1991}. There may also exist multiple stable solutions~\cite{korsak1972,turitsyn2014}. Further, non-convexities in the power flow feasible space related to multiple solutions are associated with non-zero relaxation gaps for convex relaxations of optimal power flow problems~\cite{bukhsh_tps,pscc2014}.

The literature details a variety of methods for calculating multiple power flow solutions. These include using a range of initializations for a Newton-based method~\mbox{\cite{overbye1996,overbye2000}}, a semidefinite relaxation~\cite{allerton2011,lavaei_convexpf}, an auxiliary gradient approach~\cite{guedes2005}, monotone operator theory~\cite{dj2015}, and a holomorphic embedding of the power flow equations~\cite{feng2015}. While these approaches can often find multiple power flow solutions, they are not guaranteed to find \emph{all} solutions. 

The numerical continuation method of~\cite{thorp1993} 
scales with the \emph{actual} number of power flow solutions
rather than an upper bound on the \emph{potential} number of solutions
meaning that it is computationally tractable for large test cases. 
Although it is claimed in~\cite{thorp1993} that this method will find 
all power flow solutions for all systems, \cite{chen2011} shows that the 
proof of this claim is flawed with an explicit counterexample presented in~\cite{counterexample2013}.  
A recent modification of this method~\cite{lesieutre_wu_allerton2015} formulates the power flow equations as intersecting ellipsoids which results in all continuation traces being bounded. The method in~\cite{lesieutre_wu_allerton2015} finds all solutions to a variety of small systems, including the counterexample in~\cite{counterexample2013}.
However, there is currently no known proof showing that this method always 
succeeds in finding all power flow solutions.

There are several methods which are guaranteed to find all power flow solutions. These methods formulate the power flow equations as polynomials whose variables are the real and imaginary parts of the voltage phasors. This enables the application of algorithms which find all \emph{complex} solutions to these polynomials. Relevant algorithms include interval analysis~\cite{mori1999}, Gr\"obner bases~\cite{ning2009}, an eigenvalue technique~\cite{dreesen2009}, ``moment/sum-of-squares'' relaxations~\cite{lasserre_book}, and numerical polynomial homotopy continuation (NPHC)~\cite{salam1989,guo1990,mehta2014a,chandra2015exploring,chandra2015equilibria}. 

While these methods find all complex solutions, only the real solutions are physically meaningful. Since the number of complex solutions is typically much greater than the number of real solutions, these methods are computationally limited to small systems. Of the numerical methods known to find all power flow solutions, NPHC is currently the most computationally tractable. NPHC has been applied to power system test cases with up to 14 buses~\cite{mehta2014a}, and up to the equivalent of 18 buses for the related Kuramoto model~\cite{mehta2015algebraic}. 

NPHC computes solutions to the power flow equations 
by tracking so-called solution paths.
The number of such paths is based on an upper bound on the number of
complex solutions to the power flow equations.  
Thus, this method can be improved by deriving tighter upper bounds.
Moreover, such tighter upper bounds improve the characterization of power flow 
feasible spaces and are therefore both interesting in their own right and 
relevant to, e.g., analyzing convex relaxations~\cite{pscc2014,bukhsh_tps}. 
For an $n$-bus system, one upper bound on the number of complex solutions
to the power flow equations is based on B\'ezout's theorem,
namely $2^{2n-2}$ complex solutions~\cite{baillieul1982}. By exploiting structure
in the equations of lossless systems of PV buses, 
this upper bound can be reduced to $\binom{2n - 2}{n-1} = \frac{(2n-2)!}{((n-1)!)^2}$
\cite{baillieul1982}. This bound was extended to (potentially lossy) systems of PQ buses in~\cite{li1987numerical} and to general power systems in~\cite{guo1994}. (See~\cite{marecek2014} for an alternative proof of this bound.)

Whereas the bounds in~\cite{baillieul1982,li1987numerical,guo1994,marecek2014} are network agnostic, the monomial structure of the power flow equations is determined by the network topology. Topology-dependent upper bounds have the potential to be tighter than previous bounds for specific problems. Recent work~\cite{acc2016} uses extensive numerical experiments to verify the topology-dependent bound in~\cite{guo1990} and conjecture a new bound for a different class of topologies. Related work~\cite{chen2015bounds} proposes topology-dependent bounds for a broader class of network topologies. For generic parameter values, the approach in~\cite{chen2015bounds} tightly bounds the number of complex solutions.

All known existing approaches have studied the maximum number of complex solutions in order to bound the number of real solutions. This paper addresses the question of whether systems exist for which the number of real solutions is equal to the 
number of complex solutions.\footnote{This was first posed as Question~5.1 in~\cite{baillieul1982}.} Otherwise, the equations may 
have additional structure that could be exploited 
to produce tighter bounds on the maximum number of real solutions.
Two-bus systems can have two real solutions and there exist three-bus systems with six real solutions\mbox{\cite{tavora1972,baillieul1984,klos1991}}, both of which match the upper bound provided by the number of complex solutions. Thus, this question is answered in the affirmative for these systems. For a lossless, four-bus system of PV buses with unity voltage magnitudes and a restricted set of line susceptances, the maximum number of real solutions is 14~\cite{baillieul1982}, which is strictly less than the bound of~20 provided 
by the number~of~complex~solutions. 

We consider lossless, four-bus systems of PV buses connected by lines with arbitrary susceptances. The first contribution of this paper is a choice of susceptance parameters that has 16 real solutions, thus demonstrating that open regions in the parameter space 
exist with more than 14 real solutions.
This choice, and many others, arose by solving 100,000 random instances
of the susceptance parameters for systems with zero power injections.
The second contribution of this paper is 
an analysis of the Galois group which allows for one to count
the number of real solutions of the power flow equations in this particular case,
namely lossless, four-bus systems of PV buses with unity voltage magnitues
and zero power injections, by simply counting the number of positive 
roots of a univariate sextic polynomial.
Although we are unable to explicitly compute this polynomial where the coefficients
are functions of the parameters, we can compute it explicitly 
when the parameters are specified.  
This analysis and computations lead to the conjecture 
that no choice of parameters will produce more than 
16 isolated real solutions. If this conjecture holds, this would provide 
a negative answer to the question of whether 
the number of real solutions can be equal to the bound provided by the number of 
complex solutions. In particular, this suggests that there is the potential for developing 
tighter bounds on the maximum number of real solutions to the power flow equations 
beyond the bounds provided by the number of complex solutions.

The rest of this paper is organized as follows. Section~\ref{l:powerflow} overviews the power flow equations and presents the four-bus test case studied in this paper 
along with a choice of parameters, resulting from a computational experiment, that yields 16 real solutions. Section~\ref{l:galois} reviews computing Galois groups and provides
the reduction down to counting the number of positive roots of 
a univariate polynomial.  
Section~\ref{l:conclusion} concludes the paper and discusses future work.

\section{The Power Flow Equations}
\label{l:powerflow}

This section provides an overview of the power flow equations
and then presents the lossless, four-bus system of PV buses which 
is the main focus of this paper. 
We conclude this section with a summary of a computational experiment
and a choice of parameters for which the power flow equations have 16 real solutions.

\subsection{Overview of the Power Flow Equations}
Consider an $n$-bus electric power system with buses labeled 
$1,\ldots,n$. The network topology and electrical parameters are described
via the admittance matrix $\mathbf{Y} = \mathbf{G} + \mathbf{j} \mathbf{B}$, where $\mathbf{j} = \sqrt{-1}$. (See, e.g.,~\cite{glover_sarma_overbye} for details on the admittance matrix construction.) Lossless systems have~$\mathbf{G} = 0$. 

Each bus has two associated complex quantities: the active and reactive power injections $P + \mathbf{j} Q$ where $P,Q \in\mathbb{R}^n$ and the voltage phasors $V_{d} + \mathbf{j} V_{q}$ denoted using rectangular coordinates with $V_{d},V_{q}\in\mathbb{R}^n$. Let $\left| V_i \right|$ denote the voltage magnitude at bus~$i$. The power flow equations are quadratic polynomials in the voltage components:

\vskip -0.1in
{\footnotesize
\begin{subequations}\label{power_flow_equations}
\begin{align}\label{Prect} 
& P_i = V_{di} \sum_{k=1}^n \left( \mathbf{G}_{ik} V_{dk} - \mathbf{B}_{ik} V_{qk} \right) + V_{qi} \sum_{k=1}^n \left( \mathbf{B}_{ik}V_{dk} + \mathbf{G}_{ik}V_{qk} \right) \\[-7pt] 
\label{Qrect}
& Q_i = V_{di} \sum_{k=1}^n \left( -\mathbf{B}_{ik}V_{dk} - \mathbf{G}_{ik} V_{qk}\right) \!+\! V_{qi} \sum_{k=1}^n \left( \mathbf{G}_{ik} V_{dk} - \mathbf{B}_{ik} V_{qk}\right) \\
\label{Vmag} 
& \left|V_i\right|^2 = V_{di}^2 + V_{qi}^2.
\end{align}
\end{subequations}
}
\vskip -0.15in

Each bus is classified as PQ, PV, or slack. PQ buses, which typically correspond to loads, enforce the active power~\eqref{Prect} and reactive power~\eqref{Qrect} equations with specified values for $P_i$ and $Q_i$. PV buses, which typically model generators, enforce \eqref{Prect} and \eqref{Vmag} with specified $P_i$ and $\left|V_i\right|$. The reactive power $Q_i$ may be computed as an ``output quantity'' via~\eqref{Qrect}.\footnote{This paper does not consider reactive power limited generators.} A single slack bus is selected with specified voltage magnitude $\left|V_i\right|$. The voltage phasor at the slack bus is usually chosen to provide a $0^\circ$ reference angle, i.e., $V_{di} = \left|V_i\right|$ and $ V_{qi} = 0$. The active power $P_i$ and reactive power $Q_i$ at the slack bus are determined from \eqref{Prect} and \eqref{Qrect}, thus ensuring conservation of power.

With the complex voltage phasors separated into real and imaginary parts $V_{d}$ and $V_{q}$, a solution to~\eqref{power_flow_equations} is only physically meaningful if all variables $V_{d}$ and $V_{q}$ are \emph{real valued}.

\subsection{Four-Bus Test Case}
\label{l:fourbus}

\begin{figure}[ht]
\centering
\begin{circuitikz}[scale=0.84, transform shape]
\ctikzset{bipoles/length=0.6cm}

\path[draw,line width=4pt] (0,0) -- (0,2);
\draw (-0.4,0.4) node[below right] {$1$};
\draw (-1.6,2.7) node[right] {$V_{d1} = \left|V_1\right|$};
\draw (-1.6,2.2) node[right] {$V_{q1} = 0$};
\path[draw,line width=2pt] (-0.4,1.5) -- (0,1.5);
\path[draw,line width=2pt] (-0.7,0.5) -- (0,0.5);
\path[draw,->,line width=2pt] (-0.7,0.5) -- (-0.7,-0.5);
\draw[line width=1] (-0.8,1.5) circle (0.4);

\path[draw,line width=4pt] (5,0) -- (5,2);
\draw (5,0.4) node[below right] {$2$};
\path[draw,line width=2pt] (5,1.5) -- (5.4,1.5);
\path[draw,line width=2pt] (5,0.5) -- (5.7,0.5);
\path[draw,->,line width=2pt] (5.7,0.5) -- (5.7,-0.5);
\draw[line width=1] (5.8,1.5) circle (0.4);
\draw (5.5,2.7) node[right] {$\left| V_2 \right|$};
\draw (5.5,2.2) node[right] {$P_2$};
\path[draw,->,line width=1pt] (5.6,2.2) -- (5,2.2);

\path[draw,line width=4pt] (0,-4) -- (0,-2);
\draw (-0.4,-3.6) node[below right] {$3$};
\draw (-1.3,-1.3) node[right] {$\left|V_3\right|$};
\draw (-1.3,-1.7) node[right] {$P_3$};
\path[draw,->,line width=1pt] (-0.7,-1.7) -- (0,-1.7);
\path[draw,line width=2pt] (-0.4,-2.5) -- (0,-2.5);
\path[draw,line width=2pt] (-0.7,-3.5) -- (0,-3.5);
\path[draw,->,line width=2pt] (-0.7,-3.5) -- (-0.7,-4.5);
\draw[line width=1] (-0.8,-2.5) circle (0.4);

\path[draw,line width=4pt] (5,-4) -- (5,-2);
\draw (5,-3.6) node[below right] {$4$};
\path[draw,line width=2pt] (5,-2.5) -- (5.4,-2.5);
\draw[line width=1] (5.8,-2.5) circle (0.4);
\draw (5.5,-1.3) node[right] {$\left|V_4\right|$};
\draw (5.5,-1.7) node[right] {$P_4$};
\path[draw,->,line width=1pt] (5.6,-1.7) -- (5,-1.7);
\path[draw,line width=2pt] (5,-3.5) -- (5.7,-3.5);
\path[draw,->,line width=2pt] (5.7,-3.5) -- (5.7,-4.5);

\path[draw,line width=2pt] (0,1.5) -- (5,1.5);
\draw (2.5,1.5) node[above] {$b_{12}$};

\path[draw,line width=2pt] (0,-3.5) -- (5,-3.5);
\draw (2.5,-3.5) node[below] {$b_{34}$};

\path[draw,line width=2pt] (0,0.5) -- (1,0.5) -- (1,-2.5) -- (0,-2.5);
\draw (1,-1) node[rotate=90,anchor=south] {$b_{13}$};

\path[draw,line width=2pt] (5,0.5) -- (4,0.5) -- (4,-2.5) -- (5,-2.5);
\draw (3.45,-1) node[rotate=90,anchor=north] {$b_{24}$};

\path[draw,line width=2pt] (0,1) -- (1,1) -- (4,-3) -- (5,-3);
\draw (2.05,0.55) node[rotate=-53,anchor=north] {$b_{14}$};

\path[draw,line width=2pt] (0,-3) -- (1,-3) -- (2.3,-1.35);
\draw[line width=2pt] (2.3,-1.35) arc (-135:53:0.3);
\path[draw,line width=2pt] (2.7,-0.98) -- (4,1) -- (5,1);
\draw (3,0.5) node[rotate=--56,anchor=north] {$b_{23}$};

\end{circuitikz}
\caption{Four-Bus System}
\label{f:fourbussystem}
\end{figure}
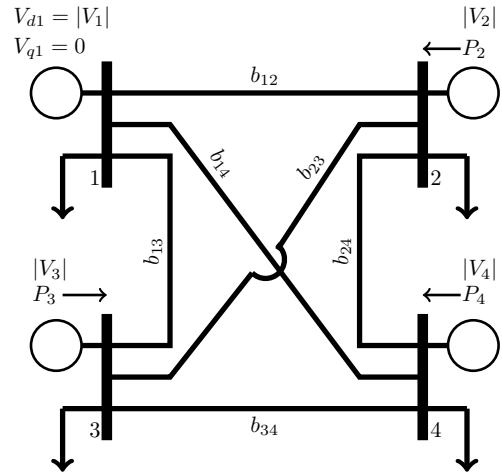

The case of interest in this paper is a lossless, four-bus system of PV buses connected by lines with arbitrary susceptances.\footnote{We follow a tradition of initially studying power flow related results for this restricted class of power systems (e.g., \cite{baillieul1982} uses a related, more specific class of systems to obtain initial results that were later generalized~\cite{li1987numerical,guo1994,marecek2014,acc2016,chen2015bounds}). Future work includes generalizing the results in this paper to four-bus systems that may have PQ buses as well as larger systems.} Fig.~\ref{f:fourbussystem} shows the one-line diagram for this system, where bus~$1$ is chosen as the slack bus. The parameter $b_{ik}$ is the susceptance of the line connecting buses~$i$ and~$k$.

The maximum number of real solutions for a system of PV buses is non-increasing with increasing resistances~\cite{baillieul1984}. Thus, allowing for losses is not necessary for the purposes of this paper: the maximum number of real solutions for lossy systems is bounded by the maximum number of real solutions for lossless systems of the same size.

With the potential for load and generation at each bus, the active power injection parameters $P_{2}$, $P_{3}$, and $P_{4}$ may be positive, negative, or zero. The voltage magnitude parameters $\left| V_i \right|$, $i=1,\ldots,4$, are positive. The line susceptance parameters may be positive (inductive), negative (capacitive) or zero.\footnote{Shunt susceptance parameters are not necessary since all buses are PV.} Note that $b_{ik} = b_{ki}$.

A susceptance parameter $b_{ik} = 0$ indicates that there is no connection between buses $i$ and $k$. Thus, the formulation we study allows for the possibility that the system is not completely connected.\footnote{The network connectivity is reflected in the monomial structure of the power flow equations. As discussed in~\cite{guo1990,acc2016,chen2015bounds}, the upper bound on the number of real power flow solutions provided by the maximum number of complex solutions is non-increasing with increasing network sparsity.}

The power flow equations for the system in Fig.~\ref{f:fourbussystem} are
\begin{subequations}
\label{pf_4bus}
\begin{align}
\nonumber & \sum_{\mathclap{\substack{k=2,3,4\\ k\neq i}}} b_{ik}\left( V_{dk}V_{qi} - V_{di}V_{qk} \right) + \left|V_1\right| V_{qi} b_{1i} = P_i \hspace{-60pt}\\[-15pt] & & i=2,3, 4 \\ \label{pf_4bus2}
& V_{di}^2 + V_{qi}^2 = \left|V_i\right|^2 & i=2,3,4 
\end{align}
\end{subequations}
%

While there are $12$ parameters in~\eqref{pf_4bus}, the three voltage magnitude 
parameters $\left|V_i\right|$ can be set to 1~per~unit without loss of generality for this system. This is observed by reformulating the power injection equations in polar coordinates:
\begin{align}\label{pf_4bus_polar}
& \sum_{\mathclap{\substack{k=1,\ldots,4\\ k\neq i}}} \left|V_i\right| \left|V_k\right| b_{ik} \sin\left(\theta_i - \theta_k\right) = P_i & i=2,3, 4
\end{align}
where $\theta_i = \arctan\left(\frac{V_{qi}}{V_{di}} \right)$ is the voltage angle at bus~$i$. Any valid change to the voltage magnitudes can be compensated by modifying the corresponding suceptance parameters such that the products of voltage magnitudes and electrical parameters are unchanged. In other words, for any choice of $\left| V_i \right|$ and $b_{ik}$, we can construct an equivalent system with susceptances $\hat{b}_{ik} = \left|V_i\right| \left| V_k\right| b_{ik} $ and 1~per~unit voltage magnitudes such 
that $\left|V_{i}\right| \left|V_k\right| b_{ik} = 1\cdot 1 \cdot \hat{b}_{ik}$.

In the particular case of lossless, four-bus systems with zero power injections ($P_i = 0$), unity voltage magnitudes ($|V_i| = 1)$, and a restricted set of susceptance parameters $b_{ik}$ (see~\cite[Prop. 2.14]{baillieul1982}),
it is shown in \cite{baillieul1982} that the maximum number of real power flow solutions is 14. 
In a computational experiment, we also take zero power injections
and unity voltage magnitues but allow the electrical parameters to be generic.
Using the NPHC method implemented in {\tt Bertini}~\cite{Bertini},
we find all complex solutions for 100,000 random instances of the 6 line 
susceptance parameters $b_{ik}$ which were selected from a normal distribution 
with mean~0 and standard deviation~8.
In this experiment, $750$ of the 100,000 instance had 16 real solutions.
Table~\ref{t:4busparameters} presents line susceptances 
for a test case that has 16 real solutions, which are listed in Table~\ref{t:4bussols}. (The characteristics of these solutions are discussed in Section~\ref{l:galois}.)
From this example, we can apply the implicit function theorem to guarantee 
that there is an open set of parameters for which there are 
16 real solutions to the power flow equations
for each value of the parameters in this set.
This statement also naturally extends to the case with non-zero power injections
as well as non-unity and non-equal voltage magnitudes.

\begin{table}[thb]
\centering
\caption{Susceptances (per unit) for a Case with 16 Real Solutions.}
\label{t:4busparameters}
{\scriptsize
\begin{tabular}{|c|c|c|c|c|c|}
\hline 
$b_{12}$ & $b_{13}$ & $b_{14}$ & $b_{23}$ & $b_{24}$ & $b_{34}$ \\\hline
1.612 & -4.649 &-5.472 &-7.504 & 10.05 & -13.571 \\ \hline
\end{tabular}
}
\end{table}

\begin{table}[thb]
\centering
\caption{The 16 Real Solutions for the System Described in Table~\ref{t:4busparameters}}
\label{t:4bussols}
{\scriptsize
\begin{tabular}{|c|c|c|c|c|c|c|c|c|}
\hline 
Sol. \# &
$V_{d2}$ & $V_{q2}$ & $V_{d3}$ & $V_{q3}$ & $V_{d4}$ & $V_{q4}$\\\hline
1 & 1 & 0 & 1& 0& -1& 0 \\ \hline
2  & -1 & 0 & 1 & 0 & -1 & 0 \\ \hline
3 &1  &0 &1 & 0&1 & 0\\ \hline
4 & 1 & 0 & -1 & 0 & -1 & 0\\ \hline
5 & -1 & 0 & 1 & 0 & 1 & 0 \\ \hline
6 & -1 & 0 & -1 & 0 & -1 & 0 \\ \hline
7 & -1 & 0& -1 & 0 & 1 & 0 \\ \hline
8 &1 & 0 & 1 & 0 &1 & 0 \\ \hline
9 & .30976 & -.95082 & -.82212 & -.56932 & -.97906 & .20359 \\ \hline
10 & -.88313& .46912 & .97310 & .23039  & .99834 & -.05754 \\ \hline
11 & .57067 & -.82118 &  -.61912 & .78530 & .41658 & -.90910 \\ \hline
12 & .89239 & -.45127 & .84624 & -.53281 & -.94751 & .31973 \\ \hline
13  & -.88313 & -.46912 & .97310 & -.23039 & .99834 & .05754 \\ \hline
14  & .57067  & .82118 & -.61912 & -.78530 & .41658 & .90910 \\ \hline
15 & .30975 & .95082 & -.82212 & .56932 & -.97906 & -.20359 \\ \hline
16 & .89239 & .45127 & .84624  & .53281  & -.94751  & -.31973\\ \hline
\end{tabular}
}
\end{table}

\section{Galois Groups and Structure}
\label{l:galois}

The test case in Section~\ref{l:fourbus} shows that there exist four-bus systems with at least 16 real solutions. We next use Galois groups to provide evidence for an upper bound of 16 real solutions for four-bus systems of PV buses.

One way to possibly understand the maximum 
number of real solutions 
to a parameterized system of polynomial equations 
is to consider the Galois group,
which corresponds geometrically to the monodromy group, e.g., see \cite{SottileBook,Harris1979}.
For a general set of the parameters, the implicit function
theorem guarantees that the solutions are analytic
functions of the parameters locally.  As one extends the domains of
these analytic functions, the Galois group encodes the relationships
between these analytic extensions.  For example, consider
the equation $x^2 - t = 0$ with 
solutions $x_{\pm}(t) = \pm \sqrt{t}$.  Although
these functions are not analytic at the origin, one can extend them
to analytic functions on, say, the complex unit circle $t(\theta) = e^{j\theta}$.
Hence, $x_{\pm}(t(\theta)) = \pm e^{j\theta/2}$.  These analytic
extensions are the same since $x_\pm(t(\theta)) = x_-(t(\theta+2\pi))$.
In this case, the Galois group is the full symmetric group
on $2$ elements, $S_2$, i.e., the plus and minus 
branches of the square root function are equivalent.

The size of the Galois group is inversely related to the structure
in the solutions.  That is, the smaller the Galois group, 
the more structure there is.  
The structure identified by the Galois group 
can then be exploited,
for example, to help solve the system more efficiently
or to aid in the derivation of upper bounds on the number of real solutions.

We start with the lossless, four-bus systems of PV buses with 
unity voltage magnitudes $(|V_i| = 1)$.  
With this, the system defined by \eqref{pf_4bus} has 
9 free parameters (six line susceptances $b_{ik}$ and three active power injections $P_i$) and 6~variables ($V_{di}$, $V_{qi}$, $i=2,3,4$).
The values in Table~\ref{t:4busparameters} with $P_i=0$
corresponds to an instance of the parameters $p^*\in\bC^{9}$.
The bound from~\cite{baillieul1982} indicates that there exist a maximum of $\binom{2n - 2}{n-1} = \binom{\,6\,}{\,3\,} = 20$ complex solutions for this system
which is sharp for parameters $p^*$.  Upon fixing an ordering
of the $20$ complex solutions, the relationship between the Galois group 
and the monodromy group describe all the ways that the $20$ solutions can 
be permuted, i.e.,~reordered, as the parameters are moved along a loop starting and ending at $p^*$.
This computation can be performed using numerical algebraic geometry,
e.g., see \mbox{\cite{BertiniBook,LeykinSottile}}.  
In this case, the Galois group is as large as possible, namely 
the symmetric group on~$20$ elements,~$S_{20}$.  That is,
for any permutation, i.e.., reordering, of the $20$ solutions, there is a loop
starting and ending at $p^*$ which realizes this permutation.  
As in the case with the square root function, each of the $20$
solution branches are equivalent meaning that there is no 
structure that we can ascertain from the Galois group~in~this~case.

We now restrict to the case of zero active power injections, 
i.e., $P_i = 0$ for $i = 2,3,4$. 
While systems with zero active power injections are a special case, 
the number of real power flow solutions decreases as the power injections approach the network's maximum loadability limit, beyond which no real solutions exist~\cite{overbye1996}. We therefore expect that the maximum number of real solutions 
is achieved with small or zero active power injections. 
Future work includes generalizing the following results 
to cases of active power injections that are not~near~zero.

After specifying zero active power injections, the remaining parameter
space is $6$ dimensional.  The approach in \cite{GaloisHRS}
yields that the Galois group is a subgroup of $S_{20}$ of order 46,080. 
Since $S_{20}$ has order $20! \approx 2.43 \cdot 10^{18}$, 
this Galois group is degenerate
which means that there is a significant amount of 
structure in the solutions as a function of the six parameters that 
we now aim to determine and exploit.  

There are two reasons for this degeneracy, both of which can be observed in the real solutions in Table~\ref{t:4bussols}.
First, eight of the solutions are always $V_{qi} = 0$
and \mbox{$V_{di} = \pm 1$} for \mbox{$i = 2,3,4$}. 
These solutions are associated with all choices of voltage angle differences of $0^\circ$ or $180^\circ$ between each pair of buses. This results in the sine of the angle differences in~\eqref{pf_4bus_polar} being equal to zero, and thus both zero active power injections and zero active power flows in the network.

Second, the other $12$ solutions are in a two-way symmetry,
i.e., if $(V_{d},V_{q})$ is a solution, so is 
$(V_{d},-V_{q})$. 
This arises from the odd symmetry of the sine function in \eqref{pf_4bus_polar}. With angle differences that are not equal to $0^\circ$ or $180^\circ$, active power flows in loops around the network in these solutions, but no active power is \emph{injected} at any of the buses. Related phenomena have been observed in practical power systems~\cite{janssens2003}.

These two observations impose restrictions on the Galois group.  
The first shows that only $12$ of the 20 complex solutions are nonconstant
with respect to the parameters.  The second shows that these $12$
nonconstant solutions arise in six pairs.  The largest
order that the Galois group could be with these two restrictions is 
$2^6\cdot 6! = \mbox{46,080}$.
Since the Galois group indeed has this order, these two observations
completely describe the degeneracy.

We ignore the $8$ solutions that are constant with respect to parameters
and focus on the other $12$ (potentially complex) solutions.
Since every coordinate $V_{qi}$ for $i = 2,3,4$ is generically distinct
for these $12$ solutions, we select the coordinate $V_{q4}$.
With this setup, there is a univariate sextic polynomial $r_6(x)$ 
whose coefficients are polynomials of the parameters $b_{ik}$ 
such that $r_6(V_{q4}^2) = 0$ exactly on the $12$ solutions.  That is,
the $V_{q4}$ coordinates of the $12$ solutions arise by
taking the positive and negative square roots of the solutions of $r_6(x) = 0$.
Therefore, the real solutions of the original system can be counted and computed
from the positive roots of $r_6(x)$.
For example, with the setup from
Table~\ref{t:4busparameters}, the sextic polynomial $r_6$ is
$$\begin{array}{r} 
\hbox{\small $x^6+13.4913 x^{5}+136.2685 x^4 - 144.4123 x^3$} \\
\hbox{\small $+~18.9004 x^2-0.5871 x+0.0017$}.\end{array}$$
By Descartes' rule of signs, $r_6$ has at most $4$ positive roots,
which is sharp in this case, yielding that the original system
has exactly $8+2\cdot4 = 16$ real solutions as listed in Table~\ref{t:4bussols}.

Suppose now that we additionally set one of the parameters 
$b_{ik}$ to zero, i.e., no line directly connects buses~$i$ and~$k$.
Then, the system generically has only $16$ complex solutions.\footnote{A maximum of $16$ complex solutions for this case was first shown by a similar argument in~\cite[Prop. 3.2]{baillieul1982}. We show that all $16$ solutions can be real for this case.}
For example, consider taking $b_{12} = 0$ and the other parameters 
as in Table~\ref{t:4busparameters}.
With this setup, the $8$ nontrivial solutions correspond with the square roots
of the roots of a quartic polynomial $r_4$, namely
$$\hbox{\small $x^4 - 1.438x^3 + 0.611x^2 - 0.070x + 0.002$}.$$
Since all four roots are positive, the corresponding system
has $8+2\cdot 4 = 16$ real solutions showing that
all $16$ complex solutions can all~be~real.

There are $20$ complex solutions when the parameters $b_{ik}$ are generic,
but only $16$ complex solutions when one of them is zero
and the others are generic. Therefore, for concreteness,
we consider the case as we take the parameter $b_{12}$ towards $0$.  
Since there are only $16$ solutions with $b_{12} =0$,
the other four solutions must diverge to infinity
as $b_{12}\rightarrow0$.  The solutions at infinity
are not identically zero and satisfy the system of homogeneous 
equations resulting from the highest degree terms in each of the polynomials.
For example, the finite solutions satisfy $V_{di}^2 + V_{qi}^2 = 1$
for $i = 2,3,4$, while the infinite solutions must satisfy
the corresponding homogeneous equations $V_{di}^2 + V_{qi}^2 = 0$.\footnote{
For $f(x,y) = x^2 + y^2 - 1$, the homogenization of $f$ 
is the polynomial $f^h(x,y,z) = z^2 f(x/z,y/z) = x^2 + y^2 - z^2$.  
Solutions at infinity have $z = 0$, i.e., 
satisfy $f^h(x,y,0) = x^2 + y^2 = 0$.}
Since the only real solution 
satisfying $V_{di}^2 + V_{qi}^2 = 0$ for $i = 2,3,4$ is the origin, 
the $4$ solutions at infinity arising from the four divergent paths
must be nonreal.  In particular, these four solutions must
be nonreal for all values of $b_{12}$ near zero.  The same
argument applies to any parameter $b_{ik}$ near zero.  

In summary, in the lossless, four-bus systems of PV buses with 
unity voltage magnitudes and zero power injections, 
we have shown that the maximum number of real solutions
is $16$ when any of the parameters $b_{ik}$ is at or near zero (i.e., a pair of buses $i$ and $k$ have a ``weak'' direct connection).  
The implicit function theorem extends this result to more general cases.
For example, in the lossless, four-bus systems of PV buses with 
unity voltage magnitudes, the maximum number of real solutions
is $16$ when the power injections are at or near zero and 
any one of the parameters $b_{ik}$ is at or near~zero.  

Due to this result together with facts that no real solutions exist for sufficiently large injections and that the number of real solutions do not increase with increasing resistances~\cite{baillieul1984},
we conjecture that the upper bound of $16$ real solutions extends
to general four-bus systems of PV buses. This conjecture is further supported by the numerical experiment described in Section~\ref{l:fourbus},which did not restrict any of the line susceptances to be near zero.

\section{Conclusions and Future Work}
\label{l:conclusion}

The power flow equations are key to many power system computations, and the characteristics of their solutions are relevant to problems in dynamics and optimization. The two main contributions of this paper are related to the maximum number of real power flow solutions. First, a test case with 16 real solutions surpasses the previously known maximum of 14 real solutions for a four-bus system. Second, Galois group theory is used to 
analyze the system and leads to the conjecture that the power flow equations for a lossless, four-bus system of PV buses has no more than 16 real solutions. 
If true, this conjecture suggests that the number of real solutions for this class of problems is strictly less than the upper bound of 20 provided by the maximum number of complex solutions. This result provides evidence for an affirmative answer regarding the open question of whether there exists a gap between the maximum number of real solutions and the upper bound provided by the number of complex solutions for systems with more than three buses.

Future work includes completing the proof of the conjecture that there are at most 16 real solutions for lossless, four-bus systems of PV buses by considering cases that do not restrict any of the parameters to be near zero. Other future work includes considering systems with PQ buses as well as developing tighter upper bounds on the number of real power flow solutions for larger test cases and other network topologies.

\bibliographystyle{IEEEtran}
\bibliography{cdc2016_nrealsol}{}

\end{document}